\documentclass[11pt,a4paper]{amsart}
\usepackage[T1]{fontenc}
\usepackage{lmodern,amsmath,amssymb,amsthm,mathtools,booktabs,array,enumitem,microtype}
\usepackage[hidelinks]{hyperref}
\hypersetup{pdftitle={Lefschetz properties of generalized whiskered graph algebras in arbitrary characteristic},pdfauthor={Tran Quang Hoa}}
\numberwithin{equation}{section}
\newtheorem{theorem}{Theorem}[section]
\newtheorem{lemma}[theorem]{Lemma}
\newtheorem{proposition}[theorem]{Proposition}
\newtheorem{corollary}[theorem]{Corollary}
\theoremstyle{definition}
\newtheorem{definition}[theorem]{Definition}
\newtheorem{example}[theorem]{Example}
\theoremstyle{remark}

\newcommand{\kk}{\Bbbk}
\newcommand{\B}{\mathcal B}
\newcommand{\Ind}{\operatorname{Ind}}
\newcommand{\floor}[1]{\left\lfloor #1\right\rfloor}
\newcommand{\ceil}[1]{\left\lceil #1\right\rceil}
\title[Generalized whiskered graph algebras]{Lefschetz properties of generalized whiskered graph algebras in arbitrary characteristic}
\author{Tran Quang Hoa}
\address{University of Education, Hue University, 34 Le Loi Street, Hue, Vietnam}
\email{tranquanghoa@hueuni.edu.vn}
\date{}
\begin{document}
\begin{abstract}
For a finite simple graph $G$ and a field $\kk$, let
\[
A(G)=\frac{\kk[x_v:v\in V(G)]}
{(x_v^2:v\in V(G))+(x_ux_v:uv\in E(G))}.
\]
We characterize the weak and strong Lefschetz properties of $A(G)$ over arbitrary fields when $G$ is obtained by attaching a positive number of leaves to each vertex of a graph. Our main tool is a decomposition of $A(G)$ as a graded $\kk[\ell]$-module into shifted quadratic monomial complete intersections indexed by the independent sets of the base graph, where $\ell$ is the sum of the variables. For ordinary whiskerings, this decomposition proves the conjecture of Holleben and Nicklasson that the weak Lefschetz property fails when the base graph has independence number at least three. We also characterize the weak Lefschetz property of the Nagata idealizations of ordinary whiskered graph algebras in characteristic zero.
\end{abstract}
\maketitle
\section{Introduction}
Let $G$ be a finite simple graph and $\kk$ a field. The Artinian algebra associated with $G$ is defined by
\[
A(G)=\frac{\kk[x_v:v\in V(G)]}{(x_v^2:v\in V(G))+(x_vx_w:vw\in E(G))}.
\]
A subset of $V(G)$ is independent if no two of its vertices are adjacent. The independence number $\alpha(G)$ is the maximum size of an independent set of $G$. A monomial basis of $A(G)_i$ is indexed by the independent sets of $G$ of size $i$. Hence the Hilbert series of $A(G)$ coincides with the independence polynomial of $G$,
\[
I(G;x)=\sum_{i=0}^{\alpha(G)}s_i(G)x^i,
\]
where $s_i(G)$ denotes the number of independent sets of size $i$ \cite{NguyenTran}. In this paper, we study the weak and strong Lefschetz properties (WLP and SLP) of these algebras for graphs obtained by attaching a positive number of leaves to each vertex of an arbitrary graph.

For background on the Lefschetz properties, we refer the reader to \cite{HarimaBook}. In characteristic zero, monomial complete intersections have the SLP \cite{PhuongTran,Stanley}. In positive characteristic, the situation is more subtle; see \cite{Cook,JonssonKling,LundqvistNicklasson}. The Lefschetz properties of quadratic monomial quotients and squarefree reductions have been studied in \cite{DaoNair,HollebenMixed,HollebenRees,JonssonKlingQuadratic,MiglioreNagelSchenck}. For the algebras $A(G)$, these properties have been investigated when $G$ is a path or a cycle \cite{NguyenTran}, a lollipop graph \cite{HoaPhuocSonLollipop}, or a tadpole graph \cite{HoaPhuocSonTadpole}. We also refer to \cite{GasanovaLundqvistNicklasson} for basis splittings of monomial algebras.

For a graph $H$ on $m\geq1$ vertices, let $w(H)$ denote the graph obtained by attaching one leaf to each vertex of $H$. Cooper et al.~\cite[Question~1.1 and Theorem~1.2]{Cooper} studied the WLP of $A(w(H))$ and established maximal rank for multiplication maps in certain degrees. Holleben and Nicklasson~\cite[Propositions~3.2 and~3.5]{HollebenNicklasson} proved that $A(w(H))$ has the SLP when $H$ is complete and the characteristic is zero or greater than $m$, and the WLP in characteristic zero when $\alpha(H)\leq2$ and $m$ is odd. They also established failure of surjectivity in certain degrees and failure of the WLP when $\alpha(H)\geq m/3+2$ \cite[Theorem~3.10 and Corollary~3.12]{HollebenNicklasson}. They conjectured that $A(w(H))$ fails the WLP whenever $\alpha(H)\geq3$ and verified the conjecture for bipartite graphs, partly by computation \cite[Conjecture~3.15 and Remark~3.16]{HollebenNicklasson}.

More generally, let $V(H)=\{v_1,\ldots,v_m\}$ and let $\mathbf p=(p_1,\ldots,p_m)$, where each $p_i$ is a positive integer. Let $H(\mathbf p)$ be the graph obtained from $H$ by attaching $p_i$ leaves to each vertex $v_i$, and set $t=\sum_{i=1}^m p_i$. Let $\Ind(H)$ denote the independence complex of $H$, whose faces are the independent sets of $H$. For $q\geq0$, we set
\[
\B_q=\kk[z_1,\ldots,z_q]/(z_1^2,\ldots,z_q^2),
\qquad \B_0=\kk.
\]
Throughout the paper, $\ell$ denotes the sum of the variables in the algebra under consideration. We regard these algebras as graded $\kk[\ell]$-modules, with $\ell$ acting by multiplication. For a graded module $M$, we use the convention $M(-s)_i=M_{i-s}$.

Our main result (Theorem~\ref{thm:attachments}) gives an isomorphism of graded $\kk[\ell]$-modules over any field,
\[
A(H(\mathbf p))\cong
\bigoplus_{S\in\Ind(H)}
\B_{t-\sum_{v_i\in S}p_i}(-|S|).
\]
For ordinary whiskerings, this refines the block structure described in \cite[Lemma~3.1]{HollebenNicklasson} for complete base graphs.

The shifted Boolean summands give explicit WLP and SLP criteria in terms of the leaf multiplicities and the independent sets of $H$ (Corollaries~\ref{cor:attachments} and~\ref{cor:structural}), with the additional conditions $p>\ceil{t/2}$ and $p>t$, respectively, in characteristic $p>0$. They also yield rank formulas for all powers of $\ell$ (Corollary~\ref{cor:attachranks}). A weighted independence enumerator determines the graded Jordan degree type of multiplication by $\ell$ over any fixed field, and the converse holds in characteristic zero or in characteristic greater than $t$ (Proposition~\ref{prop:weightedJDT}).

For ordinary whiskerings, two summands suffice to prove the conjecture of Holleben and Nicklasson over any field (Theorem~\ref{thm:conjecture}). Theorem~\ref{thm:whiskered} determines exactly which multiplication maps have maximal rank. As further applications, we classify connected graphs $G$ of girth at least six for which $A(G)$ is level and has the WLP, and show that, in characteristic zero, the Nagata idealization of $A(w(H))$ has the WLP if and only if $A(w(H))$ does (Theorems~\ref{thm:levelgirth} and~\ref{thm:idealization}).

The paper is organized as follows. Section~\ref{sec:boolean} develops the Boolean module criteria. Section~\ref{sec:whisker} treats generalized whiskerings, and Section~\ref{sec:ordinary} specializes to ordinary whiskerings. Section~\ref{sec:applications} gives applications to level and Gorenstein algebras.

\section{Boolean modules and Lefschetz criteria}\label{sec:boolean}
Let $P_r$ and $C_r$ denote the path and cycle on $r$ vertices, respectively, with $r\geq3$ for cycles. For the graph $\varnothing$ with no vertices, we set $\alpha(\varnothing)=0$ and $A(\varnothing)=\kk$. For $q\geq0$, we adopt the convention that $\binom qr=0$ whenever $r<0$ or $r>q$.

\begin{definition}
A standard graded Artinian $\kk$-algebra $A$ has the weak Lefschetz property (WLP) if there exists a linear form $\ell\in A_1$ such that the multiplication map $\ell:A_i\to A_{i+1}$ has maximal rank for every $i\geq 0$. It has the strong Lefschetz property (SLP) if there exists a linear form $\ell\in A_1$ such that the multiplication map $\ell^{j-i}:A_i\to A_j$ has maximal rank for every $0\leq i<j$. Such a linear form is called a weak or strong Lefschetz element, respectively.
\end{definition}
For Artinian algebras defined by monomial ideals, it suffices to consider the sum of the variables when testing the Lefschetz properties.

\begin{theorem}[\cite{LundqvistNicklasson,MiglioreMiroRoigNagel,Nicklasson}]\label{thm:canonical}
	Let $R=\kk[x_1,\ldots,x_n]$ be a polynomial ring over a field $\kk$, let $I\subset R$ be an Artinian monomial ideal, and set $\ell=x_1+\cdots+x_n$. Then $R/I$ has the WLP if and only if $\ell$ is a weak Lefschetz element. Similarly, $R/I$ has the SLP if and only if $\ell$ is a strong Lefschetz element.
\end{theorem}
In fact, under the hypotheses of Theorem~\ref{thm:canonical}, the canonical linear form $\ell$ maximizes the rank of each multiplication map $\ell^{j-i}:(R/I)_i\to(R/I)_j$ separately. To see this, we may extend $\kk$ to an infinite field. For fixed degrees $i<j$, the linear forms for which multiplication by their $(j-i)$-th power attains the largest possible rank form a nonempty Zariski open subset of $R_1$. This subset contains a linear form with all coefficients nonzero. Since $I$ is a monomial ideal, rescaling the variables induces a graded automorphism of $R/I$ that takes this linear form to $\ell$. As field extensions preserve ranks, the assertion follows.

For $A(G)$ and $\B_q$, we use the canonical linear forms
\[
\ell=\sum_{v\in V(G)}x_v\quad\text{and}\quad
\ell=\sum_{i=1}^q z_i,
\]
respectively. We regard these algebras as graded $\kk[\ell]$-modules, where the indeterminate $\ell$ acts by multiplication by the corresponding canonical linear form.

For a finite-dimensional graded $\kk[\ell]$-module $M$, a homogeneous Jordan string with initial degree $a$ and final degree $b$ corresponds to the graded cyclic module $(\kk[\ell]/(\ell^{b-a+1}))(-a)$. The multiset of intervals $[a,b]$ arising from a homogeneous Jordan decomposition of $M$ is called its \emph{graded Jordan degree type}; see \cite{AltafiIarrobinoMaciasMarques}. It is uniquely determined by the ranks of the maps $\ell^{j-i}:M_i\to M_j$ for all $i\leq j$. In characteristic zero, by the SLP of $\B_q$ \cite{PhuongTran,Stanley}, its Jordan strings correspond to the intervals $[i,q-i]$ with multiplicity $\binom{q}{i}-\binom{q}{i-1}$ for $0\leq i\leq\floor{q/2}$. We call $a+b$ the \emph{doubled center} of the interval $[a,b]$.

\begin{theorem}\label{thm:criterion}
	Suppose that $\kk$ has characteristic zero and that a standard graded Artinian $\kk$-algebra $A$ admits a decomposition
	\[
	A\cong\bigoplus_\lambda\B_{q_\lambda}(-s_\lambda)
	\]
	as graded $\kk[\ell]$-modules for a linear form $\ell\in A_1$ that tests both Lefschetz properties. Set
	\[
	c_-=\min_\lambda(q_\lambda+2s_\lambda),\qquad
	c_+=\max_\lambda(q_\lambda+2s_\lambda).
	\]
	Then
	\begin{align}
		A\text{ has the WLP}&\iff \floor{c_+/2}\leq\ceil{c_-/2},\label{eq:wlp}\\
		A\text{ has the SLP}&\iff c_+-c_-\leq 1.\label{eq:slp}
	\end{align}
	Equivalently, $A$ has the WLP if and only if either $c_+-c_-\leq 1$, or $c_+-c_-=2$ and $c_-$ is odd.
\end{theorem}
\begin{proof}
Let $\ell$ be the linear form in the hypothesis. Let $\sigma$ and $\epsilon$ be the largest initial degree and the smallest final degree, respectively, of the homogeneous Jordan strings of $A$. We claim that $\ell$ is a weak Lefschetz element if and only if $\sigma\leq\epsilon$. Indeed, if $\sigma\leq\epsilon$, then multiplication by $\ell$ is injective in degrees $i<\epsilon$ and surjective in degrees $i\geq\epsilon$. Conversely, if $\ell$ is a weak Lefschetz element, then the map $\ell:A_\epsilon\to A_{\epsilon+1}$ has nonzero kernel and hence is surjective. Since $A/(\ell)$ is standard graded, surjectivity propagates to all higher degrees. Thus no Jordan string starts after degree $\epsilon$, giving $\sigma\leq\epsilon$.

Similarly, $\ell$ is a strong Lefschetz element if and only if the Jordan intervals are pairwise nested. Indeed, two nonnested intervals $[a,b]$ and $[c,d]$ with $a<c$ and $b<d$ give rise to a nonzero kernel and a nonzero cokernel of the map $\ell^{d-a}:A_a\to A_d$. Conversely, if a map $\ell^{j-i}:A_i\to A_j$ with $i<j$ has both a nonzero kernel and a nonzero cokernel, then there exist two Jordan intervals, one containing $i$ but not $j$ and the other containing $j$ but not $i$. These intervals cannot be nested.

In characteristic zero, the shortest Jordan interval in $\B_q(-s)$ is
\[
[\floor{(q+2s)/2},\ceil{(q+2s)/2}],
\]
and all other Jordan intervals in this summand contain it. Taking the largest initial degree and the smallest final degree over all summands, we obtain
\[
\sigma=\floor{c_+/2},\qquad \epsilon=\ceil{c_-/2},
\]
which proves \eqref{eq:wlp}.

If $c_+-c_-\leq 1$, the doubled centers of any two Jordan intervals differ by at most one. Since their endpoints are integers, two nonnested intervals would have doubled centers differing by at least two. Thus the Jordan intervals are pairwise nested. Conversely, if $c_+-c_-\geq 2$, the shortest intervals in the summands with doubled centers $c_-$ and $c_+$ have strictly increasing initial and final degrees, and hence are not nested. This proves \eqref{eq:slp}. The last assertion follows from \eqref{eq:wlp} by considering the parity of $c_-$.
\end{proof}

For $q\geq 0$ and integers $u,v$ with $0\leq u\leq v\leq q$, let $R_q(u,v)$ denote the rank of the multiplication map $\ell^{v-u}:(\B_q)_u\to(\B_q)_v$, and set $R_q(u,v)=0$ otherwise. For $0\leq u\leq v\leq q$, define
\[
(a,d)=
\begin{cases}
	(u,v),&\text{if }u+v\leq q,\\
	(q-v,q-u),&\text{if }u+v>q.
\end{cases}
\]
With respect to the monomial bases, the multiplication matrix is $(v-u)!$ times the inclusion matrix between subsets of sizes $u$ and $v$. By taking complements and transposing if necessary, we reduce to the inclusion matrix between subsets of sizes $a$ and $d$, where $a+d\leq q$. Wilson's integral diagonal form for inclusion matrices \cite{Wilson} gives
\begin{equation}\label{eq:boolrank}
	R_q(u,v)=\sum_{h=0}^{a}\left(\binom{q}{h}-\binom{q}{h-1}\right)
	\mathbf{1}_{\kk}\left((v-u)!\binom{d-h}{a-h}\right).
\end{equation}
Here, for an integer $n$, $\mathbf{1}_{\kk}(n)=1$ if $n$ is nonzero in $\kk$, and $\mathbf{1}_{\kk}(n)=0$ otherwise. Since Wilson's diagonal form is obtained by unimodular row and column transformations over $\mathbb Z$, the formula is valid over any field.

In characteristic $p>0$, Jonsson Kling \cite[Theorem~2.9]{JonssonKling} established the following maximal rank criterion for $0\leq u<v\leq q$:
\begin{equation}\label{eq:powerrank}
	R_q(u,v)=\min\left\{\binom{q}{u},\binom{q}{v}\right\}
	\quad\Longleftrightarrow\quad p>\min\{v,q-u\}.
\end{equation}
For completeness, we give a short proof using \eqref{eq:boolrank}. Set $r=v-u=d-a>0$. Since $a\leq\floor{q/2}$, every multiplicity in \eqref{eq:boolrank} is positive, and their sum is $\binom{q}{a}=\min\{\binom{q}{u},\binom{q}{v}\}$. The corresponding diagonal entries are
\[
r!\binom{d-h}{a-h}=\prod_{b=0}^{r-1}(d-h-b).
\]
If $p>d$, none of these entries vanishes in $\kk$. If $p\leq r$, all these entries vanish in $\kk$ since $r!=0$ in $\kk$. If $r<p\leq d$, take $h=d-p$. Then $0\leq h\leq a$, and the corresponding diagonal entry has a factor equal to $p$. Since $d=\min\{v,q-u\}$, this proves \eqref{eq:powerrank}. In particular, for $0\leq i<q$, the multiplication map $\ell:(\B_q)_i\to(\B_q)_{i+1}$ has maximal rank if and only if
\begin{equation}\label{eq:adjrank}
	p>\min\{i+1,q-i\}.
\end{equation}

The following lemma collects the known characteristic bounds for the WLP and SLP of $\B_q$; see \cite{Cook,JonssonKling,LundqvistNicklasson,PhuongTran}.

\begin{lemma}\label{lem:boolpositive}
Let $q\geq 0$ and suppose that $\operatorname{char}\kk=p>0$. Then $\B_q$ has the WLP if and only if $p>\ceil{q/2}$, and has the SLP if and only if $p>q$. Moreover, if $p>\ceil{q/2}$, all multiplication maps by $\ell$ have the same ranks as in characteristic zero. If $p>q$, the same holds for multiplication by every power of $\ell$.
\end{lemma}

\begin{proof}
	For $q=0$, the assertions are immediate. By Theorem~\ref{thm:canonical}, it suffices to consider the canonical linear form $\ell$. For $q\geq 1$, the WLP assertion follows from \eqref{eq:adjrank}, since $\max_{0\leq i<q}\min\{i+1,q-i\}=\ceil{q/2}$. The SLP assertion follows from \eqref{eq:powerrank}, since $\max_{0\leq u<v\leq q}\min\{v,q-u\}=q$. In the respective characteristic ranges, all the relevant multiplication maps have maximal rank, as in characteristic zero.
\end{proof}

\section{Lefschetz properties of generalized whiskered graphs}\label{sec:whisker}
Let $H$ be a finite simple graph with vertex set $V(H)=\{v_1,\ldots,v_m\}$, where $m\geq 1$. For $\mathbf p=(p_1,\ldots,p_m)$ with $p_i\geq 1$, let $H(\mathbf p)$ be the graph obtained from $H$ by attaching $p_i$ leaves to each vertex $v_i$, and set $t=\sum_{i=1}^m p_i$.

\begin{theorem}\label{thm:attachments}
	Over an arbitrary field $\kk$, there is an isomorphism of graded $\kk[\ell]$-modules
	\begin{equation}\label{eq:attach}
		A(H(\mathbf p))\cong
		\bigoplus_{S\in\Ind(H)}\B_{t-\sum_{v_i\in S}p_i}(-|S|).
	\end{equation}
\end{theorem}
\begin{proof}
Write $x_i=x_{v_i}$ and let $y_{i,1},\ldots,y_{i,p_i}$ be the variables corresponding to the leaves attached to $v_i$. Consider the invertible change of variables
\[
z_{i,1}=y_{i,1}+x_i,\qquad z_{i,j}=y_{i,j}\quad (2\leq j\leq p_i),
\]
with the variables $x_i$ unchanged. Since $x_i^2=y_{i,1}^2=x_iy_{i,1}=0$, we have $z_{i,1}^2=x_iz_{i,1}=0$. Conversely, the original relations follow by substituting $y_{i,1}=z_{i,1}-x_i$. All other defining relations remain unchanged. Thus, in the new variables, the defining relations are
\[
x_i^2=0,\quad x_ix_h=0\ \text{if }v_iv_h\in E(H),\quad z_{i,j}^2=x_iz_{i,j}=0.
\]
Moreover, the canonical linear form becomes
\[
\ell=\sum_{i=1}^m\sum_{j=1}^{p_i}z_{i,j}.
\]
We partition the standard monomial basis in the new variables according to the set $S$ of vertices $v_i$ for which $x_i$ occurs in a monomial. By the defining relations, each such set $S$ is independent in $H$. For each $S\in\Ind(H)$, let $M_S$ be the subspace spanned by the monomials
\[
\left(\prod_{v_i\in S}x_i\right)
\left(\prod_{v_i\notin S}\prod_{j\in J_i}z_{i,j}\right),
\qquad J_i\subseteq\{1,\ldots,p_i\}.
\]

Since $\ell$ involves only the variables $z_{i,j}$, the subspace $M_S$ is invariant under multiplication by $\ell$. Multiplication by $\prod_{v_i\in S}x_i$ identifies the Boolean algebra generated by the variables $z_{i,j}$ for $v_i\notin S$ with $M_S$, up to a grading shift by $|S|$. This identification is compatible with multiplication by $\ell$. Hence
\[
M_S\cong\B_{t-\sum_{v_i\in S}p_i}(-|S|)
\]
as graded $\kk[\ell]$-modules. Taking the direct sum over all $S\in\Ind(H)$ gives \eqref{eq:attach}.
\end{proof}

As an application of Theorem~\ref{thm:attachments}, we prove the conjecture of Holleben and Nicklasson \cite[Conjecture~3.15]{HollebenNicklasson}.

\begin{theorem}\label{thm:conjecture}
Let $H$ be a finite simple graph with $\alpha(H)\geq 3$. Then $A(w(H))$ does not have the WLP over any field.
\end{theorem}
\begin{proof}
Let $m=|V(H)|$, $i=\ceil{m/2}$, and $s=2i+2-m$. Then $s=2$ if $m$ is even and $s=3$ if $m$ is odd. Since $\alpha(H)\geq 3$, there exists an independent set $S$ of size $s$. By Theorem~\ref{thm:attachments}, the decomposition of $A(w(H))$ contains the summands $\B_m$ and $\B_{m-s}(-s)$, corresponding to $\varnothing$ and $S$, respectively. Since
\[
\binom{m}{i}>\binom{m}{i+1},
\]
the multiplication map $\ell:(\B_m)_i\to(\B_m)_{i+1}$ has nonzero kernel. On the other hand, since $m+s=2i+2$, we have
\[
\binom{m-s}{i-s}<\binom{m-s}{i+1-s}.
\]
Indeed, $i+1-s=(m-s)/2$ is the middle degree of $\B_{m-s}$. Thus the multiplication map on the summand $\B_{m-s}(-s)$ from degree $i$ to degree $i+1$ has nonzero cokernel. Consequently, the multiplication map $\ell:A(w(H))_i\to A(w(H))_{i+1}$ has both nonzero kernel and nonzero cokernel, and hence fails to have maximal rank. By Theorem~\ref{thm:canonical}, $A(w(H))$ does not have the WLP.
\end{proof}
For a nonnegative weight function $w$ on $V(H)$, let
\[
\alpha_w(H)=\max_{S\in\Ind(H)}\sum_{v\in S}w(v)
\]
denote the weighted independence number of $H$. Writing $x_+=\max\{x,0\}$, we set
\[
a_+=\alpha_{(p_i-2)_+}(H),\qquad a_-=\alpha_{(2-p_i)_+}(H),
\]
where the weights are assigned to the corresponding vertices $v_i$.

\begin{corollary}\label{cor:attachments}
In characteristic zero, $A(H(\mathbf p))$ has the WLP if and only if
\[
\floor{\frac{t+a_-}{2}}\leq\ceil{\frac{t-a_+}{2}}.
\]
It has the SLP if and only if $a_++a_-\leq 1$. In characteristic $p>0$, the same criteria hold with the additional conditions $p>\ceil{t/2}$ for the WLP and $p>t$ for the SLP, respectively.
\end{corollary}
\begin{proof}
The doubled center of the summand indexed by $S\in\Ind(H)$ is
\[
t+\sum_{v_i\in S}(2-p_i).
\]
Since every subset of an independent set is independent, vertices with negative contributions can be omitted when maximizing the sum, and vertices with positive contributions can be omitted when minimizing it. Hence the minimum and maximum of these doubled centers are $t-a_+$ and $t+a_-$, respectively. The assertions in characteristic zero follow from Theorem~\ref{thm:criterion}.

Suppose that $\operatorname{char}\kk=p>0$. The summand indexed by the empty set is $\B_t$, and every other Boolean factor has at most $t$ variables. By Lemma~\ref{lem:boolpositive}, if $p\leq\ceil{t/2}$ (resp.\ $p\leq t$), some multiplication map by $\ell$ (resp.\ by a power of $\ell$) in the summand $\B_t$ has both nonzero kernel and nonzero cokernel. The corresponding map on $A(H(\mathbf p))$ therefore fails to have maximal rank. By Theorem~\ref{thm:canonical}, $A(H(\mathbf p))$ fails to have the WLP (resp.\ SLP). If $p>\ceil{t/2}$ (resp.\ $p>t$), all multiplication maps by $\ell$ (resp.\ by powers of $\ell$) in every summand have the same ranks as in characteristic zero. Thus the corresponding characteristic-zero criterion remains valid in positive characteristic.
\end{proof}

For each $j\geq 1$, let $X_j=\{v_i:p_i=j\}$. We regard the empty set as a clique.

\begin{corollary}\label{cor:structural}
	In characteristic zero, $A(H(\mathbf p))$ has the SLP if and only if one of the following conditions holds:
	\begin{enumerate}[label=\textup{(\roman*)}]
		\item $p_i=2$ for all $i$;
		\item $p_i\in\{2,3\}$ for all $i$, and $X_3$ is a nonempty clique;
		\item $p_i\in\{1,2\}$ for all $i$, and $X_1$ is a nonempty clique.
	\end{enumerate}
	The algebra $A(H(\mathbf p))$ has the WLP if and only if one of the conditions \textup{(i)--(iii)} holds or one of the following conditions holds:
	\begin{enumerate}[label=\textup{(\roman*)},start=4]
		\item $t$ is odd, $p_i\in\{2,3,4\}$ for all $i$, $\alpha(H[X_3])\leq 2$, and $X_4$ is a clique whose vertices are adjacent to every vertex of $X_3$. Moreover, either $X_4\neq\varnothing$ or $\alpha(H[X_3])=2$;
		\item $t$ is even, $p_i\in\{1,2,3\}$ for all $i$, and both $X_1$ and $X_3$ are nonempty cliques;
		\item $t$ is odd, $p_i\in\{1,2\}$ for all $i$, and $\alpha(H[X_1])=2$.
	\end{enumerate}
	In characteristic $p>0$, the SLP holds if and only if one of \textup{(i)--(iii)} holds and $p>t$, whereas the WLP holds if and only if one of \textup{(i)--(vi)} holds and $p>\ceil{t/2}$.
\end{corollary}
\begin{proof}
By definition, $a_-=\alpha(H[X_1])$. Thus $a_-=0$ if and only if $X_1=\varnothing$, and $a_-=1$ if and only if $X_1$ is a nonempty clique. Similarly, $a_+=0$ if and only if all $p_i\leq 2$, and $a_+=1$ if and only if all $p_i\leq 3$ and $X_3$ is a nonempty clique. Therefore, the pairs $(a_+,a_-)=(0,0),(1,0),(0,1)$ correspond to \textup{(i)--(iii)}, respectively, and give exactly the SLP cases in characteristic zero by Corollary~\ref{cor:attachments}.

By Theorem~\ref{thm:criterion}, the remaining WLP cases have $a_++a_-=2$ and $t-a_+$ odd. If $(a_+,a_-)=(2,0)$, all $p_i$ belong to $\{2,3,4\}$. Conversely, assuming $p_i\in\{2,3,4\}$ for all $i$, we have $a_-=0$, and the bound $a_+\leq 2$ is equivalent to $\alpha(H[X_3])\leq 2$, $X_4$ being a clique, and every vertex of $X_4$ being adjacent to every vertex of $X_3$. Under these conditions, equality holds if and only if $X_4\neq\varnothing$ or $\alpha(H[X_3])=2$. Together with odd $t$, this gives \textup{(iv)}. For $(a_+,a_-)=(1,1)$, all $p_i$ belong to $\{1,2,3\}$, both $X_1$ and $X_3$ are nonempty cliques, and $t$ is even, giving \textup{(v)}. For $(a_+,a_-)=(0,2)$, all $p_i$ belong to $\{1,2\}$, $\alpha(H[X_1])=2$, and $t$ is odd, giving \textup{(vi)}. The positive-characteristic conditions follow from Corollary~\ref{cor:attachments}.
\end{proof}
\begin{corollary}\label{cor:attachranks}
Over any field, for $0\leq i\leq j$, we have
\[
\operatorname{rank}(\ell^{j-i}:A(H(\mathbf p))_i\to A(H(\mathbf p))_j)
=\sum_{S\in\Ind(H)}R_{t-\sum_{v_h\in S}p_h}(i-|S|,j-|S|),
\]
where $R_q(u,v)$ is defined in Section~\ref{sec:boolean} and computed by \eqref{eq:boolrank} for $0\leq u\leq v\leq q$. In particular, for the ordinary whiskering $w(H)$, we have
\begin{equation}\label{eq:whiskerrank}
\operatorname{rank}(\ell^{j-i}:A(w(H))_i\to A(w(H))_j)
=\sum_{r=0}^{\alpha(H)}s_r(H)R_{m-r}(i-r,j-r).
\end{equation}
\end{corollary}

\begin{proof}
By Theorem~\ref{thm:attachments}, the rank of the multiplication map is the sum of its ranks on the individual summands, taking into account the grading shifts. For ordinary whiskerings, $p_h=1$ for all $h$ and $t=m$. Grouping the summands according to the cardinality of $S$ gives \eqref{eq:whiskerrank}.
\end{proof}

\begin{corollary}\label{cor:whiskerjdt}
Let $H$ and $H'$ be finite simple graphs. Over a fixed field $\kk$, the following conditions are equivalent:
\begin{enumerate}[label=\textup{(\roman*)}]
\item $H$ and $H'$ have the same independence polynomial;
\item $A(w(H))$ and $A(w(H'))$ are isomorphic as graded $\kk[\ell]$-modules;
\item $A(w(H))$ and $A(w(H'))$ have the same graded Jordan degree type.
\end{enumerate}
\end{corollary}
\begin{proof}
\textup{(i)} $\Rightarrow$ \textup{(ii)}. Suppose that $H$ and $H'$ have the same independence polynomial. Then $|V(H)|=|V(H')|$ and $s_r(H)=s_r(H')$ for all $r\geq 0$. If both graphs are empty, the assertion is immediate. Otherwise, by Theorem~\ref{thm:attachments}, the graded $\kk[\ell]$-modules $A(w(H))$ and $A(w(H'))$ are isomorphic.

\textup{(ii)} $\Rightarrow$ \textup{(i)}. Suppose that the graded modules $A(w(H))$ and $A(w(H'))$ are isomorphic. Then they have the same Hilbert series, so $I(w(H);x)=I(w(H');x)$. Let $m=|V(H)|=|V(H')|$, since the coefficient of $x$ in the independence polynomial of a whiskered graph is twice the order of its base graph. If $m=0$, then both graphs have no vertices and the assertion is immediate. The decomposition in Theorem~\ref{thm:attachments} gives the classical formula for the independence polynomial of a whiskered graph (see \cite{LevitMandrescu}):
\[
I(w(H);x)=(1+x)^m I\left(H;\frac{x}{1+x}\right).
\]
The Hilbert series of $A(w(H))$ therefore determines the independence polynomial
\[
I(H;x)=(1-x)^m I\left(w(H);\frac{x}{1-x}\right).
\]
Thus $H$ and $H'$ have the same independence polynomial.

\textup{(ii)} $\Leftrightarrow$ \textup{(iii)}. The graded Jordan degree type is a complete invariant of finite-dimensional graded $\kk[\ell]$-modules.
\end{proof}

For arbitrary leaf multiplicities, we define the bivariate weighted independence enumerator
\[
F_{H,\mathbf p}(x,y)=\sum_{S\in\Ind(H)}x^{|S|}y^{\sum_{v_i\in S}p_i}\in\mathbb Z[x,y].
\]

\begin{proposition}\label{prop:weightedJDT}
Let $t=\sum_{i=1}^m p_i$. Over a fixed field $\kk$, the polynomial $F_{H,\mathbf p}$ determines the graded $\kk[\ell]$-module $A(H(\mathbf p))$ up to isomorphism. Conversely, if $\kk$ has characteristic zero or characteristic $p>t$, the graded Jordan degree type of $A(H(\mathbf p))$ determines $F_{H,\mathbf p}$.
\end{proposition}
\begin{proof}
The coefficient of $x$ in $F_{H,\mathbf p}(x,y)$ is $\sum_{i=1}^m y^{p_i}$, which determines $t=\sum_{i=1}^m p_i$. Moreover, the coefficient of $x^sy^w$ is precisely the multiplicity of $\B_{t-w}(-s)$ in the decomposition of Theorem~\ref{thm:attachments}. This proves the first assertion.

Conversely, the highest nonzero degree of $A(H(\mathbf p))$ is $t$. Indeed, all $t$ attached leaves form an independent set, and each original vertex belonging to an independent set excludes at least one of its attached leaves. Thus $\alpha(H(\mathbf p))=t$, and hence the graded Jordan degree type determines $t$.

If $\operatorname{char}\kk=p>t$, Lemma~\ref{lem:boolpositive} shows that every Boolean summand has the same graded Jordan degree type as in characteristic zero. In either case, we group the Jordan intervals according to their doubled centers. For a fixed doubled center $c$, consider the summands $\B_{c-2s}(-s)$ in increasing order of $s$, where $0\leq s\leq\floor{c/2}$. Each such summand contains exactly one interval $[s,c-s]$, and all its other intervals start in degrees greater than $s$. After removing the intervals contributed by summands with smaller values of $s$, the remaining multiplicity of $[s,c-s]$ is precisely the multiplicity of $\B_{c-2s}(-s)$. Proceeding in this way, we recover the multiplicity of every Boolean summand and hence all coefficients of $F_{H,\mathbf p}$, since the multiplicity of $\B_{c-2s}(-s)$ corresponds to the monomial $x^sy^{t-c+2s}$ whenever this summand occurs. This proves the converse.
\end{proof}

The following example shows that the converse of Proposition~\ref{prop:weightedJDT} may fail in positive characteristic.

\begin{example}\label{ex:weighted-char2}
	Let $H=P_4$ with vertices in path order, and consider $\mathbf p=(4,2,1,2)$ and $\mathbf p'=(2,3,1,3)$. Both vectors have total weight $t=9$, and a direct computation gives
	\[
	F_{H,\mathbf p}(x,y)-F_{H,\mathbf p'}(x,y)
	=xy^2(1-y)(1-y-xy).
	\]
	Since this difference vanishes at $y=(1+x)^{-1}$ and the Hilbert series is $(1+x)^t F_{H,\mathbf p}(x,(1+x)^{-1})$ by Theorem~\ref{thm:attachments}, the two algebras have the same Hilbert series.
	
	Suppose that $\operatorname{char}\kk=2$. For $q\geq1$, replacing $z_1$ by $\ell=z_1+\cdots+z_q$ gives
	\[
	\B_q\cong
	\kk[\ell]/(\ell^2)\otimes_\kk\B_{q-1},
	\]
	where $\ell$ acts on the first factor. Thus all Jordan strings have length two. Every independent set of $P_4$ is proper and all leaf multiplicities are positive, so every Boolean summand has at least one variable. Thus the number of Jordan strings starting in degree $r$ is the coefficient of $x^r$ in the quotient of the Hilbert series by $1+x$. Hence the two algebras have the same graded Jordan degree type but different weighted independence enumerators.
\end{example}

\section{Ordinary whiskerings and applications}\label{sec:ordinary}

We now apply the results of the preceding section to ordinary whiskerings and determine which multiplication maps have maximal rank.
\begin{theorem}\label{thm:whiskered}
Let $H$ be a graph on $m\geq 1$ vertices and set $\alpha=\alpha(H)$. Over any field $\kk$, for $0\leq i<j\leq m$, the multiplication map
\[
\ell^{j-i}:A(w(H))_i\longrightarrow A(w(H))_j
\]
fails to have maximal rank if and only if at least one of the following conditions holds:
\begin{enumerate}[label=\textup{(\roman*)}]
\item $m<i+j<m+\alpha$;
\item $\operatorname{char}\kk=p>0$ and $p\leq\min\{j,m-i\}$.
\end{enumerate}
	
In particular, in characteristic zero, $A(w(H))$ has the WLP if and only if $\alpha=1$, or $\alpha=2$ and $m$ is odd, and has the SLP if and only if $H$ is complete. In characteristic $p>0$, the same criteria hold with the additional conditions $p>\ceil{m/2}$ and $p>m$, respectively.
\end{theorem}

\begin{proof}
By Theorem~\ref{thm:attachments}, we have
\begin{equation}\label{eq:whiskered}
	A(w(H))\cong\bigoplus_{S\in\Ind(H)}\B_{m-|S|}(-|S|).
\end{equation}
	
Suppose first that $\operatorname{char}\kk=0$. For an independent set $S$ of size $s$, the corresponding summand $\B_{m-s}(-s)$ is supported in degrees $s,\ldots,m$ and has doubled center $m+s$. By the symmetry and unimodality of its Hilbert function and the SLP of $\B_{m-s}$, the multiplication map on this summand is injective if $i+j\leq m+s$ and surjective if $i+j\geq m+s$. These assertions remain valid when $i<s$, since in this case the source is zero and $i+j<m+s$. Consequently, the multiplication map on $A(w(H))$ is injective if $i+j\leq m$ and surjective if $i+j\geq m+\alpha$.
	
Now suppose that $m<i+j<m+\alpha$. Since $\binom mi>\binom mj$, the multiplication map on the summand $\B_m$ has nonzero kernel. Set $s=i+j+1-m$. Then $2\leq s\leq\alpha$ and $s\leq j$. Choose an independent set $S$ of size $s$ and consider the corresponding summand $\B_{m-s}(-s)$. If $i<s$, the source is zero and the target is nonzero. If $i\geq s$, then
\[
\binom{m-s}{i-s}<\binom{m-s}{j-s},
\]
since $i<j$ and $i+j<m+s$. In either case, the multiplication map on this summand has nonzero cokernel. Thus the multiplication map on $A(w(H))$ has both nonzero kernel and nonzero cokernel. Since this argument uses only the dimensions of the graded components, condition \textup{(i)} implies failure of maximal rank in every characteristic.
	
Suppose now that $\operatorname{char}\kk=p>0$. If condition \textup{(ii)} holds, then by \eqref{eq:powerrank}, the multiplication map on the summand $\B_m$ has both nonzero kernel and nonzero cokernel. Hence the multiplication map on $A(w(H))$ fails to have maximal rank. Otherwise, we have $p>\min\{j,m-i\}$. For a summand $\B_{m-s}(-s)$ with both degrees in its support, the characteristic bound in \eqref{eq:powerrank} is
\[
\min\{j-s,m-i\}\leq\min\{j,m-i\}.
\]
It follows that the multiplication map on each summand has the same rank as in characteristic zero. This is also true when either the source or the target is zero. Consequently, the multiplication map on $A(w(H))$ has the same rank as in characteristic zero. This proves the maximal rank criterion.
	
Finally, in the notation of Corollary~\ref{cor:attachments}, we have $t=m$, $a_+=0$, and $a_-=\alpha$. Thus, in characteristic zero, the WLP criterion is
\[
\floor{(m+\alpha)/2}\leq\ceil{m/2},
\]
which is equivalent to $\alpha=1$, or $\alpha=2$ and $m$ odd. The SLP criterion is $\alpha\leq1$, which is equivalent to $H$ being complete. The additional characteristic bounds follow from the same corollary.
\end{proof}

Theorem~\ref{thm:whiskered} completely answers Question~1.1 of Cooper et al.~\cite{Cooper} for ordinary whiskered graphs over any field, including the case where the base graph has no edges.

We next consider the case where exactly two leaves are attached to each vertex of the base graph.
\begin{corollary}\label{cor:twoleaves}
	Let $H^{[2]}$ be the graph obtained from a graph $H$ on $m\geq1$ vertices by attaching two leaves to each vertex. Then $A(H^{[2]})$ has the SLP if and only if $\operatorname{char}\kk=0$ or $\operatorname{char}\kk=p>2m$. Moreover, its Hilbert function is symmetric about degree $m$ in every characteristic.
\end{corollary}
\begin{proof}
Since $p_i=2$ for all $i$, we have $t=2m$ and $a_+=a_-=0$. The SLP criterion follows from Corollary~\ref{cor:attachments}. By Theorem~\ref{thm:attachments}, each summand $\B_{2m-2|S|}(-|S|)$ has a Hilbert function symmetric about degree $m$, and so does their direct sum.
\end{proof}
As a further application, we obtain trees of every order whose associated algebras have the SLP.

\begin{corollary}\label{cor:orders}
In characteristic zero, for every positive integer $n$, there exists a tree $T$ on $n$ vertices with maximum degree at most four such that $A(T)$ has the SLP.
\end{corollary}

\begin{proof}
	For $n=1$, take $T=P_1$, whose associated algebra has the SLP. For $n\geq2$, write $n=3m+j$, where $m\geq1$ and $j\in\{-1,0,1\}$. Let $v$ be an endpoint of $P_m$ (the unique vertex if $m=1$). Attach two leaves to each vertex of $P_m$ other than $v$, and attach $2+j$ leaves to $v$. The resulting tree $T$ has $3m+j=n$ vertices and maximum degree at most four. The total number of attached leaves is $t=2m+j$.
	
	By Theorem~\ref{thm:attachments}, the doubled center of the summand indexed by $S\in\Ind(P_m)$ is
	\[
	t-\sum_{v_i\in S}(p_i-2)
	=
	\begin{cases}
		2m & \text{if } v\in S,\\
		2m+j & \text{if } v\notin S.
	\end{cases}
	\]
Thus $c_+-c_-=|j|\leq1$, and Theorem~\ref{thm:criterion} implies that $A(T)$ has the SLP.
\end{proof}

\section{Applications to level and Gorenstein algebras}\label{sec:applications}

The socle of $A(G)$ has a monomial basis indexed by the maximal independent sets of $G$. Thus $A(G)$ is level if and only if $G$ is well-covered, that is, all maximal independent sets of $G$ have the same cardinality. The girth of a graph is the length of its shortest cycle, with the convention that forests have infinite girth.

\begin{theorem}\label{thm:levelgirth}
Let $G$ be a connected graph of girth at least six such that $A(G)$ is level. Then $A(G)$ has the WLP if and only if one of the following conditions holds:
\begin{enumerate}[label=\textup{(\roman*)}]
\item $G\cong K_1$, $K_2$, or $P_4$;
\item $G\cong C_7$ or $w(P_3)$, and $\operatorname{char}\kk\neq2$.
\end{enumerate}
\end{theorem}
\begin{proof}
Since $A(G)$ is level, $G$ is well-covered. By the classification of Finbow, Hartnell, and Nowakowski~\cite{FinbowHartnellNowakowski}, either $G\cong K_1$ or $C_7$, or $G$ has a perfect matching consisting of pendant edges. In the latter case, deleting one leaf from each matched pair gives a connected graph $H$ of girth at least six such that $G\cong w(H)$.
	
Suppose that $A(w(H))$ has the WLP. By Theorem~\ref{thm:whiskered}, either $\alpha(H)=1$, or $\alpha(H)=2$ and $|V(H)|$ is odd. In the first case, $H$ is complete. Since its girth is at least six, it contains no triangles, and hence $H\cong K_1$ or $K_2$. In the second case, if $H$ has at least five vertices, then it contains a connected induced subgraph on five vertices. Since the girth of $H$ is at least six, this subgraph is a tree. Being bipartite, it has an independent set of size at least three, contradicting $\alpha(H)=2$. Since $|V(H)|$ is odd and $\alpha(H)=2$, we conclude that $H$ has exactly three vertices. As $H$ is connected and contains no triangles, we have $H\cong P_3$.
	
Conversely, by Theorem~\ref{thm:whiskered}, $A(w(K_1))$ and $A(w(K_2))$ have the WLP over any field, whereas $A(w(P_3))$ has the WLP if and only if $\operatorname{char}\kk\neq2$. Note that $w(K_1)\cong K_2$ and $w(K_2)\cong P_4$.
	
For $G=K_1$, the assertion is immediate. For $G=C_7$, label the vertices by $\mathbb Z/7\mathbb Z$. If $\sum_i c_ix_i$ belongs to the kernel of $\ell:A(C_7)_1\to A(C_7)_2$, then $c_i+c_{i+2}=0$ for every $i$. Iterating these relations seven times gives $c_i=-c_i$. Thus the map is injective if $\operatorname{char}\kk\neq2$. In characteristic two, $\ell^2=0$ and $\ell\neq0$, so the map is not injective. Since the source and target have dimensions $7$ and $14$, respectively, the map has maximal rank if and only if $\operatorname{char}\kk\neq2$.
	
Every independent set of size three in $C_7$ is of the form $\{i,i+2,i+4\}$. Since
	\[
	\ell x_ix_{i+4}=x_ix_{i+2}x_{i+4},
	\]
the multiplication map $\ell:A(C_7)_2\to A(C_7)_3$ is surjective over any field. All remaining multiplication maps between consecutive degrees have maximal rank. Therefore, by Theorem~\ref{thm:canonical}, $A(C_7)$ has the WLP if and only if $\operatorname{char}\kk\neq2$.
\end{proof}

From now on, we assume that $\operatorname{char}\kk=0$. Let $A$ be a standard graded Artinian level $\kk$-algebra of socle degree $m$, and let $A^\vee$ denote its graded dual, with $(A^\vee)_j=(A_{-j})^*$. Its Nagata idealization is
\[
B=A\ltimes A^\vee(-m-1),
\]
with graded components $B_i=A_i\oplus(A_{m+1-i})^*$ and multiplication $(a,f)(b,g)=(ab,ag+bf)$, where $(ag)(c)=g(ac)$. Since $A$ is level, $A^\vee(-m-1)$ is generated as an $A$-module in degree one, so $B$ is standard graded. Moreover, $B$ is Artinian Gorenstein of socle degree $m+1$.

Let $H$ be a graph on $m$ vertices. Every maximal independent set of $w(H)$ contains exactly one vertex from each whisker pair. Hence $A(w(H))$ is level of socle degree $m$. Its Nagata idealization is G-quadratic and hence Koszul by the results of D'Al\`i and Venturello~\cite{DaliVenturello}; see also \cite[Theorem~4.4 and Lemma~4.5]{HollebenNicklasson}.

Dao and Nair~\cite{DaoNair} use the failure of surjectivity of certain multiplication maps in level algebras to construct Artinian Gorenstein algebras without the WLP. The following theorem extends \cite[Theorem~4.6]{HollebenNicklasson} by giving a complete characterization of the WLP for the Nagata idealizations of ordinary whiskered graph algebras.

\begin{theorem}\label{thm:idealization}
	Suppose that $\operatorname{char}\kk=0$. Let $H$ be a graph on $m\geq1$ vertices and set $\alpha=\alpha(H)$. Let $A=A(w(H))$ and let $B$ be its Nagata idealization. Then the following conditions are equivalent:
	\begin{enumerate}[label=\textup{(\roman*)}]
		\item $B$ has the WLP;
		\item $A$ has the WLP;
		\item either $\alpha=1$, or $\alpha=2$ and $m$ is odd.
	\end{enumerate}
\end{theorem}

\begin{proof}
	\textup{(ii)} $\Leftrightarrow$ \textup{(iii)}. This follows from Theorem~\ref{thm:whiskered}.
	
	\textup{(i)} $\Rightarrow$ \textup{(iii)}. Let $L\in B_1$ be a weak Lefschetz element and set $i=\ceil{m/2}$. If the map $L:B_i\to B_{i+1}$ is not surjective, it is injective by the WLP, and its dual $L:B_{m-i}\to B_{m-i+1}$ is surjective. Since $m-i\leq i$ and $B$ is standard graded, it follows that $L:B_i\to B_{i+1}$ is surjective. This is a contradiction, so the map is surjective.
	
	Let $\overline L$ be the image of $L$ under the natural graded surjection $B\to A$. Then the map $\overline L:A_i\to A_{i+1}$ is surjective. If $m$ is even and $\alpha\geq2$, or if $m$ is odd and $\alpha\geq3$, then
	\[
	m<2i+1<m+\alpha.
	\]
	By Theorem~\ref{thm:whiskered}, the map $\ell:A_i\to A_{i+1}$ fails to have maximal rank and hence is not surjective. This contradicts Theorem~\ref{thm:canonical}, since $\ell$ maximizes the rank of this map.
	
	\textup{(iii)} $\Rightarrow$ \textup{(i)}. Suppose first that $m=2k+1$ and $\alpha\leq2$. By \eqref{eq:whiskered}, the doubled centers of the summands lie between $m$ and $m+2$. Hence the maps
	\[
	\ell:A_k\to A_{k+1}
	\quad\text{and}\quad
	\ell:A_{k+1}\to A_{k+2}
	\]
	are injective and surjective, respectively. The map $(\ell,0):B_k\to B_{k+1}$ is the direct sum of the first map and the dual of the second, and hence is injective. By Gorenstein duality, the map $(\ell,0):B_{k+1}\to B_{k+2}$ is surjective. Since surjectivity propagates to higher degrees in a standard graded algebra, Gorenstein duality gives injectivity in the lower degrees. Hence $(\ell,0)$ is a weak Lefschetz element of $B$.
	
	Now suppose that $m=2k$ and $H=K_m$. By \eqref{eq:whiskered}, the map $\ell:A_k\to A_{k+1}$ is surjective, and its kernel $K$ is contained in the unshifted summand $\B_m$, since the corresponding maps on the remaining summands $\B_{m-1}(-1)$ are isomorphisms.
	
	In the coordinates introduced in the proof of Theorem~\ref{thm:attachments}, $\B_m$ is the subalgebra generated by $z_1,\ldots,z_m$. With respect to the perfect multiplication pairing
	\[
	(\B_m)_k\times(\B_m)_k
	\longrightarrow(\B_m)_m\cong\kk,
	\]
	we have $K^\perp=\ell(\B_m)_{k-1}$. Since $\B_m$ has the SLP, the map
	\[
	\ell^2:(\B_m)_{k-1}\longrightarrow(\B_m)_{k+1}
	\]
	is an isomorphism. Thus $K\cap K^\perp=0$, so the restriction of the pairing to $K$ is nondegenerate.
	
	Let $\varphi\in A_m^*$ be the functional taking the coefficient of $z_1\cdots z_m$. Suppose that $(u,f)\in B_k$ belongs to the kernel of multiplication by $(\ell,\varphi)$. Then
	\[
	\ell u=0,\qquad \ell f+u\varphi=0.
	\]
	In particular, $u\in K$. Evaluating the second equality at any $v\in K$ gives
	\[
	\varphi(uv)=-f(\ell v)=0.
	\]
	By the nondegeneracy of the pairing on $K$, we obtain $u=0$. Hence $\ell f=0$, and the surjectivity of $\ell:A_k\to A_{k+1}$ gives $f=0$. Thus the middle multiplication map $(\ell,\varphi):B_k\to B_{k+1}$ is injective and, by Gorenstein symmetry, an isomorphism. By propagation of surjectivity and Gorenstein duality, $(\ell,\varphi)$ is a weak Lefschetz element of $B$.
\end{proof}

\section*{Declaration of AI assistance}
The author used ChatGPT (OpenAI) as an auxiliary tool in developing and checking mathematical arguments and refining the exposition. The author takes full responsibility for the content of the paper.

\end{document}